\documentstyle[twoside]{article}
\title{On two Thomae-type transformations for hypergeometric series with integral parameter differences }
\author{
Y. S. Kim,\footnote{Department of Mathematics Education, Wonkwang University, Iksan, Korea
E-mail: yspkim@wonkwang.ac.kr}
\ \ Arjun. K. Rathie\footnote{Department of Mathematics, Central University of Kerala,
Kasaragad 671123, Kerala, India.
E-Mail: akrathie@cukerala.edu.in} \ \ and 
\ R. B. Paris\footnote{School of Computing, Engineering and Applied Mathematics, University of Abertay Dundee, Dundee DD1 1HG, UK.
E-Mail: r.paris@abertay.ac.uk}\ \footnote{Corresponding author}
 \\}
\begin{document}
\def\f#1#2{\mbox{${\textstyle \frac{#1}{#2}}$}}
\def\dfrac#1#2{\displaystyle{\frac{#1}{#2}}}
\def\boldal{\mbox{\boldmath $\alpha$}}
\newcommand{\bee}{\begin{equation}}
\newcommand{\ee}{\end{equation}}
\newcommand{\lam}{\lambda}
\newcommand{\ka}{\kappa}
\newcommand{\al}{\alpha}
\newcommand{\th}{\theta}
\newcommand{\om}{\omega}
\newcommand{\Om}{\Omega}
\newcommand{\fr}{\frac{1}{2}}
\newcommand{\fs}{\f{1}{2}}
\newcommand{\g}{\Gamma}
\newcommand{\br}{\biggr}
\newcommand{\bl}{\biggl}
\newcommand{\ra}{\rightarrow}
\newcommand{\mbint}{\frac{1}{2\pi i}\int_{c-\infty i}^{c+\infty i}}
\newcommand{\mbcint}{\frac{1}{2\pi i}\int_C}
\newcommand{\mboint}{\frac{1}{2\pi i}\int_{-\infty i}^{\infty i}}
\newcommand{\gtwid}{\raisebox{-.8ex}{\mbox{$\stackrel{\textstyle >}{\sim}$}}}
\newcommand{\ltwid}{\raisebox{-.8ex}{\mbox{$\stackrel{\textstyle <}{\sim}$}}}
\renewcommand{\topfraction}{0.9}
\renewcommand{\bottomfraction}{0.9}
\renewcommand{\textfraction}{0.05}
\newcommand{\mcol}{\multicolumn}
\date{}
\maketitle
\begin{abstract}
We obtain two new Thomae-type transformations for hypergeometric series with $r$ pairs of numeratorial and denominatorial parameters differing by positive integers. This is achieved by application of the so-called Beta integral method developed by Krattenthaler and Rao [Symposium on {\it Symmetries in Science} (ed. B. Gruber), Kluwer (2004)] to two recently obtained Euler-type transformations. Some special cases are given.

\vspace{0.4cm}

\noindent {\bf Mathematics Subject Classification:} 33C15, 33C20 
\vspace{0.3cm}

\noindent {\bf Keywords:} Generalized hypergeometric series, Thomae transformations, generalized Euler-type transformations 
\end{abstract}

\vspace{0.3cm}

\begin{center}
{\bf 1. \  Introduction}
\end{center}
\setcounter{section}{1}
\setcounter{equation}{0}
\renewcommand{\theequation}{\arabic{section}.\arabic{equation}}
The generalized hypergeometric function ${}_pF_q(x)$ is defined for complex parameters and argument by the series
\bee\label{e11}
{}_pF_q\left[\!\!\begin{array}{c} a_1, a_2, \ldots ,a_p\\b_1, b_2, \ldots ,b_q\end{array}\!; x\right]=
\sum_{k=0}^\infty \frac{(a_1)_k (a_2)_k \ldots (a_p)_k}{(b_1)_k (b_2)_k \ldots (b_q)_k}\,\frac{x^k}{k!}.
\ee
When $q\geq p$ this series converges for $|x|<\infty$, but when $q=p-1$ convergence occurs when $|x|<1$ (unless the series terminates).
In (\ref{e11}) the Pochhammer symbol or ascending factorial $(a)_n$ is given for integer $n$ by 
\[(a)_n=\frac{\g(a+n)}{\g(a)}=\left\{\begin{array}{ll} 1 & (n=0)\\a(a+1)\ldots (a+n-1) & (n\geq 1),\end{array}\right.\]
where $\g$ is the gamma function. In what follows we shall adopt the convention of writing the finite sequence of parameters $(a_1, a_2, \ldots ,a_p)$ simply by $(a_p)$ and the product of $p$ Pochhammer symbols by
\[((a_p))_k\equiv (a_1)_k \ldots (a_p)_k,\]
where an empty product $p=0$ is interpreted as unity.

Recent work has been carried out on the extension of various summations theorems, such as those of Gauss, Kummer, Bailey and Watson \cite{KRR, RR11, RP}, and also of Euler-type transformations to higher-order hypergeometric functions with $r$ pairs of numeratorial and denominatorial parameters differing by positive integers \cite{MP1, MP2}. Our interest in this note is concerned with obtaining  similar extensions of the two-term Thomae transformation \cite[p.~52]{S}
\[{}_3F_2\left[\!\!\begin{array}{c}a,\,b,\,c\\d,\,e\end{array}\!;1\right]=\frac{\g(d) \g(e) \g(\sigma)}{\g(a) \g(b+\sigma) \g(c+\sigma)}\,{}_3F_2\left[\!\!\begin{array}{c}c-a,\,d-a,\,\sigma\\b+\sigma,\,c+\sigma\end{array}\!;1\right]\]
for $\Re (\sigma)>0$, $\Re (a)>0$, where $\sigma=e+d-a-b-c$ is the parametric excess. 
Many other results of the above type, including three-term Thomae transformations, are given in \cite[pp.~116-121]{S}; see also \cite{W}.

The so-called Beta integral method introduced by Krathenthaler and Rao \cite{KR2} generates new identities for hypergeometric series for some fixed value of the argument (usually 1) from
known identities for hypergeometric series with a smaller number of parameters involving the argument $x$, $1-x$ or a combination of their powers. The basic idea of this method is to multiply the known hypergeometric identity by the factor $x^{d-1}(1-x)^{e-d-1}$, where $e$ and $d$ are suitable parameters, integrate term by term over $[0,1]$ making use of the beta integral representation
\bee\label{e12}
\hspace{2cm}\int_0^1 t^{a-1}(1-t)^{b-1}dt=\frac{\g(a)\g(b)}{\g(a+b)}\qquad (\Re (a)>0,
\ \Re (b)>0)
\ee
and finally to rewrite the result in terms of a new hypergeometric series. We apply this method to two Euler-type transformations  recently obtained in \cite{MP1, MP2} to derive two two-term Thomae-type transformations for hypergeometric functions with $r$ pairs of numeratorial and denominatorial parameters differing by positive integers.
\vspace{0.6cm}

\begin{center}
{\bf 2. \  Extended Thomae-type transformations}
\end{center}
\setcounter{section}{2}
\setcounter{equation}{0}
\renewcommand{\theequation}{\arabic{section}.\arabic{equation}}
Our starting point is the following Euler-type transformations for hypergeometric functions with $r$ pairs of numeratorial and denominatorial parameters differing by positive integers $(m_r)$.
\vspace{0.4cm}

\noindent{\bf Theorem 1.}\  {\it Let $(m_r)$ be a sequence of positive 
integers with $m:=m_1+\cdots +m_r$. Then
we have the two Euler-type transformations \cite{MP1, MP2} for $|\arg\,(1-x)|<\pi$
\bee\label{e21}
{}_{r+2}F_{r+1}\left[\!\!\begin{array}{c} a, b,\\c,\end{array}\!\!\!\!\begin{array}{c} (f_r+m_r)\\ (f_r)\end{array}\!; x\right]=(1-x)^{-a}\,{}_{m+2}F_{m+1}
\left[\!\!\begin{array}{c} a,\,c\!-\!b\!-\!m,\\c,\end{array}\!\!\!\!\begin{array}{c} (\xi_m+1)\\(\xi_m)\end{array}\!; \frac{x}{x-1}\right]
\ee
provided $b\neq f_j$ $(1\leq j\leq r)$, $(c-b-m)_m\neq 0$ and
\bee\label{e22}
{}_{r+2}F_{r+1}\left[\!\!\begin{array}{c} a, b,\\c,\end{array}\!\!\!\!\begin{array}{c} (f_r+m_r)\\ (f_r)\end{array}\!; x\right]
=(1-x)^{c\!-\!a\!-\!b\!-\!m}\,{}_{m+2}F_{m+1}
\left[\!\!\begin{array}{c} c\!-\!a\!-\!m,\,c\!-\!b\!-\!m,\\c,\end{array}\!\!\!\!\begin{array}{c} (\eta_m+1)\\(\eta_m)\end{array}\!; x\right]
\ee
provided $(c-a-m)_m\neq 0$, $(c-b-m)_m\neq 0$. 
The $(\xi_m)$ and $(\eta_m)$ are respectively the nonvanishing zeros of the associated parametric polynomials $Q_m(t)$ and ${\hat Q}_m(t)$ defined below.}
\vspace{0.4cm}

The parametric polynomials $Q_m(t)$ and ${\hat Q}_m(t)$, both of degree $m=m_1+\cdots +m_r$, are given by
\bee\label{e23}
Q_m(t)=\frac{1}{(\lambda)_m}\sum_{k=0}^m(b)_kC_{k,r} (t)_k (\lambda-t)_{m-k},
\ee
where $\lambda:=b-a-m$, and
\bee\label{e24}
{\hat Q}_m(t)=\sum_{k=0}^m \frac{(-1)^k C_{k,r} (a)_k (b)_k  (t)_k}{(c-a-m)_k (c-b-m)_k}\,G_{m,k}(t)
\ee
where
\[G_{m,k}(t):={}_3F_2\left[\!\!\begin{array}{c}-m+k, t+k, c-a-b-m\\c-a-m+k, c-b-m+k\end{array}\!;1\right].\]
The coefficients $C_{k,r}$ are defined for $0\leq k\leq m$ by 
\bee\label{e200}
C_{k,r}= \frac{1}{\Lambda}\sum_{j=k}^m \sigma_{j}{\bf S}_j^{(k)},\qquad \Lambda=(f_1)_{m_1}\ldots (f_r)_{m_r},
\ee
with $C_{0,r}=1$, $C_{m,r}=1/\Lambda$. The ${\bf S}_j^{(k)}$ denote the Stirling numbers of the second kind and the $\sigma_j$ $(0\leq j\leq m)$ are generated by the relation
\bee\label{e201}
(f_1+x)_{m_1} \cdots (f_r+x)_{m_r}=\sum_{j=0}^m \sigma_{j}x^j.
\ee 
For $0\leq k\leq m$, the function $G_{m,k}(t)$ is a polynomial in $t$ of degree $m-k$
and both $Q_m(t)$ and ${\hat Q}_m(t)$ are normalized so that $Q_m(0)={\hat Q}_m(0)=1$. 
\vspace{0.4cm}

\noindent{\bf Remark 1.}\ \ In \cite{MP3}, an alternative representation for the coefficients $C_{k,r}$ is given as the terminating hypergeometric series of unit argument
\[
C_{k,r}=\frac{(-1)^k}{k!}\,{}_{r+1}F_r\left[\!\!\begin{array}{c}-k,\\{}\end{array}
\!\!\!\begin{array}{c}(f_r+m_r)\\(f_r)\end{array}\!;1\right].
\]
When $r=1$, with $f_1=f$, $m_1=m$, Vandermonde's summation theorem \cite[p.~243]{S} can be used to show that
\bee\label{e203}
C_{k,1}=\left(\!\!\!\begin{array}{c}m\\k\end{array}\!\!\!\right)\,\frac{1}{(f)_k}.
\ee
\vspace{0.4cm}

We first apply the Beta integral method \cite{KR2} to the result in 
(\ref{e22}) to obtain a new hypergeometric identity.
Multiplying both sides by $x^{d-1}(1-x)^{e-d-1}$, where $e$, $d$ are arbitrary parameters satisfying $\Re (e-d)>0$, $\Re (d)>0$, we integrate over the interval $[0, 1]$. The left-hand side yields
\begin{eqnarray}
\int_0^1x^{d-1}(1-x)^{e-d-1}&&\hspace{-0.9cm} {}_{r+2}F_{r+1}\left[\!\!\begin{array}{c} a,\, b,\\c,\end{array}\!\!\!\!\begin{array}{c} (f_r+m_r)\\ (f_r)\end{array}\!; x\right]dx\nonumber\\
&=&\!\!\sum_{k=0}^\infty\frac{(a)_k (b)_k}{(c)_k\,k!}\,\frac{((f_r+m_r))_k}{((f_r))_k}\int_0^1x^{d+k-1}(1-x)^{e-d-1}dx\nonumber\\
&=&\!\!\sum_{k=0}^\infty\frac{(a)_k (b)_k}{(c)_k\,k!}\,\frac{((f_r+m_r))_k}{((f_r))_k}\,
\frac{\g(d+k) \g(e-d)}{\g(e+k)}\nonumber\\
&=&\!\!\frac{\g(d) \g(e-d)}{\g(e)} {}_{r+3}F_{r+2}\left[\!\!\begin{array}{c} a,\,b,\,d,\\c,\,e,\end{array}\!\!\!\!\begin{array}{c} (f_r+m_r)\\ (f_r)\end{array}\!; 1\right],\label{e63}
\end{eqnarray}
upon evaluation of the integral by (\ref{e12}) and use of the definition (\ref{e11}) when it is supposed that $\Re (s)>0$, where $s$ is the parametric excess given by 
\bee\label{e250}
s:=c+e-a-b-d-m.
\ee

Proceeding in a similar manner with the right-hand side of (\ref{e22}), we obtain
\begin{eqnarray}
\int_0^1&&\hspace{-0.9cm}x^{d-1}(1-x)^{s-1} {}_{m+2}F_{m+1}\left[\!\!\begin{array}{c} c-a-m,\,c-b-m,\\c,\end{array}\!\!\!\!\begin{array}{c} (\eta_m+1)\\ (\eta_m)\end{array}\!;x\right]dx\nonumber\\
&=&\!\!\sum_{k=0}^\infty\frac{(c-a-m)_k (c-b-m)_k}{(c)_k\,k!}\,\frac{((\eta_m+1))_k}{((\eta_m))_k}\int_0^1x^{d+k-1}(1-x)^{s-1}dx\nonumber\\
&=&\!\!\frac{\g(d) \g(s)}{\g(c+e-a-b-m)}
\,{}_{m+3}F_{m+2}\left[\!\!\begin{array}{c} c-a-m,\,c-b-m,\,d,\\c,\,c+e-a-b-m,\end{array}\!\!\!\!\begin{array}{c} (\eta_m+1)\\ (\eta_m)\end{array}\!; 1\right].\label{e64}
\end{eqnarray}
Then by (\ref{e63}) and (\ref{e64}) we obtain the two-term Thomae-type hypergeometric identity given in the following theorem, where the restriction $\Re (d)>0$ can be removed by appeal to analytic continuation:
\vspace{0.4cm}

\noindent{\bf Theorem 2.}\  {\it Let $(m_r)$ be a sequence of positive 
integers with $m:=m_1+\cdots +m_r$. Then
\[{}_{r+3}F_{r+2}\left[\!\!\begin{array}{c} a,\,b,\,d,\\c,\,e,\end{array}\!\!\!\!\begin{array}{c} (f_r+m_r)\\ (f_r)\end{array}\!; 1\right]\hspace{9cm}\]
\bee\label{e25}
=\frac{\g(e) \g(s)}{\g(e-d) \g(s+d)}
\,{}_{m+3}F_{m+2}\left[\!\!\begin{array}{c} c-a-m,\,c-b-m,\,d,\\c,\,s+d,\end{array}\!\!\!\!\begin{array}{c} (\eta_m+1)\\ (\eta_m)\end{array}\!; 1\right]
\ee
provided $(c-a-m)_m\neq 0$, $(c-b-m)_m\neq 0$, $\Re (e-d)>0$ and $\Re (s)>0$, where $s$ is defined by (\ref{e250}).}
\vspace{0.4cm}

The same procedure can be applied to (\ref{e21}) when the parameter $a=-n$ (to ensure convergence of the resulting integral at $x=1$), where $n$ is a non-negative integer, to yield the right-hand side of (\ref{e21}) given by
\begin{eqnarray}
\int_0^1x^{d-1}&&\hspace{-0.9cm}(1-x)^{e-d+n-1} {}_{m+2}F_{m+1}\left[\!\!\begin{array}{c} -n,\,c-b-m,\\c,\end{array}\!\!\!\!\begin{array}{c} (\xi_m+1)\\ (\xi_m)\end{array}\!;\frac{x}{x-1}\right]dx\nonumber\\
&=&\!\!\sum_{k=0}^n \frac{(-1)^k(-n)_k (c-b-m)_k}{(c)_k\,k!}\,\frac{((\xi_m+1))_k}{((\xi_m))_k} \int_0^1x^{d+k-1}(1-x)^{e-d+n-k-1}dx\nonumber\\
&=&\!\!\frac{\g(d) \g(e-d+n)}{\g(e+n)}\,\sum_{k=0}^n\frac{(-n)_k (c-b-m)_k (d)_k}{(c)_k (1-e+d-n)_k k!}\,\frac{((\xi_m+1))_k}{((\xi_m))_k}\nonumber\\
&=&\!\!\frac{\g(d) \g(e-d+n)}{\g(e+n)}\,
{}_{m+3}F_{m+2}\left[\!\!\begin{array}{c} -n,\,c-b-m,\,d,\\c,\,1-e+a+d,\end{array}\!\!\!\!\begin{array}{c} (\xi_m+1)\\ (\xi_m)\end{array}\!; 1\right]\label{e66}
\end{eqnarray}
provided $\Re (e-d)>0$, $\Re (d)>0$. From (\ref{e63}) and (\ref{e66}), and appeal to analytic continuation to remove the restriction $\Re (d)>0$, we then obtain the finite Thomae-type transformation
\vspace{0.4cm}

\noindent{\bf Theorem 3.}\  {\it Let $(m_r)$ be a sequence of positive 
integers with $m:=m_1+\cdots +m_r$. Then, for non-negative integer $n$
\bee\label{e26}
{}_{r+3}F_{r+2}\left[\!\!\begin{array}{c} -n,\,b,\,d,\\c,\,e,\end{array}\!\!\!\!\begin{array}{c} (f_r+m_r)\\ (f_r)\end{array}\!; 1\right]
=\frac{(e-d)_n}{(e)_n}
\,{}_{m+3}F_{m+2}\left[\!\!\begin{array}{c} -n,\,c-b-m,\,d,\\c,\,1-e+d-n,\end{array}\!\!\!\!\begin{array}{c} (\xi_m+1)\\ (\xi_m)\end{array}\!; 1\right]
\ee
provided $b\neq f_j$ $(1\leq j\leq r)$, $(c-b-m)_m\neq 0$ and $\Re (e-d)>0$.}
\vspace{0.6cm}

\begin{center}
{\bf 3. \  Examples}
\end{center}
\setcounter{section}{3}
\setcounter{equation}{0}
\renewcommand{\theequation}{\arabic{section}.\arabic{equation}}
When $r=0$ (with $m=0$) we recover from (\ref{e25}) and (\ref{e26}) the known results \cite{W}
\[{}_{3}F_{2}\left[\!\!\begin{array}{c} a,\,b,\,d,\\c,\,e,\end{array}\!; 1\right]=
\frac{\g(e) \g(c+e-a-b-d)}{\g(e-d) \g(c+e-a-b)}
\,{}_{3}F_{2}\left[\!\!\begin{array}{c} c-a,\,c-b,\,d\\c,\,c+e-a-b\end{array}\!; 1\right]\]
for $\Re (e-d)>0$, $\Re (e+c-a-b-d)>0$ and
\[{}_{3}F_{2}\left[\!\!\begin{array}{c} -n,\,b,\,d,\\c,\,e,\end{array}\!; 1\right]=
\frac{(e-d)_n}{(e)_n}
\,{}_{3}F_{2}\left[\!\!\begin{array}{c} -n,\,c-b,\,d,\\c,\,1-e+d-n,\end{array}\!; 1\right]\]
for $\Re (e-d)>0$ with $n$ a non-negative integer.

In the particular case $r=1$, $m_1=m=1$, $f_1=f$, we have the parametric polynomial from (\ref{e23})
\[Q_1(t)=1+\frac{(b-f)t}{(c-b-1)f}\]
with the nonvanishing zero $\xi_1=\xi$ (provided $b\neq f$, $c-b-1\neq 0$) given by
\bee\label{e30}
\xi=\frac{(c-b-1)f}{f-b},
\ee
and from (\ref{e24})
\[{\hat Q}_1(t)=1-\frac{\{(c-a-b-1)f+ab\}t}{(c-a-1)(c-b-1)f}\]
with the nonvanishing zero $\eta_1=\eta$ (provided $c-a-1\neq 0$, $c-b-1\neq 0$) given by
\bee\label{e31}
\eta=\frac{(c-a-1)(c-b-1)f}{ab+(c-a-b-1)f}.
\ee
Then we have from (\ref{e25}) and (\ref{e26}) the transformations
\[{}_{4}F_{3}\left[\!\!\begin{array}{c} a,\,b,\,d,\\c,\,e,\end{array}\!\!\!\!\begin{array}{c} f+1\\ f\end{array}\!; 1\right]
=\frac{\g(e) \g(s)}{\g(e-d) \g(s+d)}
\,{}_{4}F_{3}\left[\!\!\begin{array}{c} c-a-1,\,c-b-1,\,d,\\c,\,s+d,\end{array}\!\!\!\!\begin{array}{c} \eta+1\\ \eta\end{array}\!; 1\right]
\]
provided $c-a-1\neq 0$, $c-b-1\neq 0$, $\Re (e-d)>0$ and $\Re (s)>0$, where $s$ is defined by (\ref{e250}) with $m=1$, and
\[{}_{4}F_{3}\left[\!\!\begin{array}{c} -n,\,b,\,d,\\c,\,e,\end{array}\!\!\!\!\begin{array}{c} f+1\\ f\end{array}\!; 1\right]
=\frac{(e-d)_n}{(e)_n}
\,{}_{4}F_{3}\left[\!\!\begin{array}{c} -n,\,c-b-1,\,d,\\c,\,1-e+d-n,\end{array}\!\!\!\!\begin{array}{c} \xi+1\\ \xi\end{array}\!; 1\right]\]
for non-negative integer $n$ and $\Re (e-d)>0$.

In the case $r=1$, $m_1=2$, $f_1=f$, we have $C_{0,r}=1$, $C_{1,r}=2/f$   and $C_{2,r}=1/(f)_2$ by (\ref{e203}).
From (\ref{e23})  and (\ref{e24}) we obtain  after a little algebra the quadratic parametric polynomials $Q_2(t)$  (with zeros $\xi_1$ and $\xi_2$) and ${\hat Q}_2(t)$ (with zeros $\eta_1$ and $\eta_2$) given by 
\[Q_2(t)=1-\frac{2(f-b)t}{(c-b-2) f}+\frac{(f-b)_2t(t+1)}{(c-b-2)_2 (f)_2}\]
and
\[{\hat Q}_2(t)=1-\frac{2Bt}{(c-a-2)(c-b-2)}+\frac{Ct(1+t)}{(c-a-2)_2(c-b-2)_2},\]
where
\[ B:=\sigma'+\frac{ab}{f},\qquad C:=\sigma'(\sigma'+1)+\frac{2ab\sigma'}{f}+\frac{(a)_2(b)_2}{(f)_2}, \qquad \sigma':=c-a-b-2.\]
For example, if $a=\f{1}{4}$, $b=\f{5}{2}$, $c=\f{3}{2}$ and $f=\fs$ we have
\[Q_2(t)=1-\f{8}{3}t+\f{4}{9}t(1+t),\qquad {\hat Q}_2(t)=1+\f{16}{9}t-\f{68}{27}t(1+t),\]
whence $\xi_1=\fs$, $\xi_2=\f{9}{2}$ and $\eta_1=\fs$, $\eta_2=-\f{27}{34}$. The transformations in (\ref{e25}) and (\ref{e26}) then yield
\bee\label{e32}
{}_4F_3\left[\!\!\begin{array}{c}\vspace{0.1cm}

\f{1}{4},\,\f{5}{2},\,d,\,\f{5}{2}\\ \f{3}{2},\,e,\,\fs\end{array}\!;1\right]=\frac{\g(e) \g(e-d-\f{13}{4})}{\g(e-d) \g(e-\f{13}{4})}\,{}_4F_3\left[\!\!\begin{array}{c}\vspace{0.1cm}

-\f{3}{4},\,-3,\,d,\ \f{7}{34}\\ e-\f{13}{4},\,\fs,\,-\f{27}{34}\end{array}\!;1\right]
\ee
provided $\Re (e-d)>\f{13}{4}$, and
\bee\label{e33}
{}_4F_3\left[\!\!\begin{array}{c}\vspace{0.1cm}

-n,\,\f{5}{2},\,d,\,\f{5}{2}\\ \f{3}{2},\,e,\,\fs\end{array}\!;1\right]=\frac{(e-d)_n}{(e)_n}\,{}_4F_3\left[\!\!\begin{array}{c}\vspace{0.1cm}

-n,\,-3,\,d,\ \f{11}{2}\\ 1-e+d-n,\,\fs,\,\f{9}{2}\end{array}\!;1\right]
\ee
for non-negative integer $n$.
We remark that a contraction of the order of the hypergeometric functions on the right-hand sides of (\ref{e32}) and (\ref{e33}) has been possible since $c=\xi_1+1=\eta_1+1=\f{3}{2}$. In addition, both series on the right-hand sides terminate: the first with summation index $k=3$ and the second with index $k=\min \{n, 3\}$. A final point to mention is that for real parameters $a$, $b$, $c$ and $f$ it is possible (when $m\geq 2$) to have complex zeros.

\vspace{0.6cm}

\begin{center}
\noindent{\bf 4. \ Concluding remarks}
\end{center}
\setcounter{section}{4}
\setcounter{equation}{0}
\renewcommand{\theequation}{\arabic{section}.\arabic{equation}}
We have employed the Beta Integral method of Krattenthaler and Rao \cite{KR2} applied to two
recently obtained Euler-type transformations for hypergeometric functions with $r$ pairs of numeratorial and denominatorial parameters differing by positive integers $(m_r)$.
By this means, we have established two Thomae-type transformations given in Theorems 2 and 3.

In order to write the hypergeometric series in (\ref{e25}) and (\ref{e26}) we require the zeros 
$(\eta_m)$ and $(\xi_m)$ of the parametric polynomials ${\hat Q}_m(t)$ and $Q_m(t)$ respectively.
However, to evaluate the series on the right-hand sides of (\ref{e25}) and (\ref{e26}), {\it it is not necessary to evaluate these zeros}. This observation can be understood by reference to the hypergeometric series
\[F\equiv{}_{m+2}F_{m+1}\left[\!\!\begin{array}{c} \alpha,\,\beta,\\ \gamma,\end{array}\!\!\!\!\begin{array}{c}(\xi_m+1)\\(\xi_m)\end{array}\!;1\right]=\sum_{k=0}^\infty\frac{(\alpha)_k (\beta)_k}{(\gamma)_k\,k!}\left(1+\frac{k}{\xi_1}\right) \ldots \left(1+\frac{k}{\xi_m}\right)\]
upon use of the fact that $(a+1)_k/(a)_k=1+(k/a)$. Since the parametric polynomial $Q_m(t)$ in (\ref{e23}) can be written as $Q_m(t)=\prod_{r=1}^m \{1-(t/\xi_r)\}$ it follows that
\[F=\sum_{k=0}^\infty\frac{(\alpha)_k (\beta)_k}{(\gamma)_k\,k!}\,Q_m(-k).\]
Consequently it is sufficient to know only the parametric polynomial $Q_m(t)$.
A similar remark applies to the series involving the zeros $(\eta_m)$ with the parametric polynomial $Q_m(-k)$ replaced by ${\hat Q}_m(-k)$.
\vspace{0.6cm}

\noindent{\bf Acknowledgement:}\ \ \ Y. S. Kim acknowledges the support of the Wonkwang University Research Fund (2013).


\vspace{0.6cm}

\begin{thebibliography}{99}
\footnotesize{

\bibitem{KRR}
Y. S. Kim, M. A. Rakha and A. K. Rathie, Extensions of certain classical summation theorems for the series ${}_2F_1$, ${}_3F_2$ and ${}_4F_3$ with applications in Ramanujan summations, Int. J. Math. Math. Sci.  309503, 26 pages (2010).

\bibitem{KR2}
C. Krattenthaler and K. Srinivasa Rao, On group theoretical aspects, hypergeometric transformations and symmetries of angular momentum coefficients, XIII Symposium on {\it Symmetries in Science}, B. Gruber (ed.), Kluwer Academic Publishers, 2004.

\bibitem{MP1}
A. R. Miller and R. B. Paris, Euler-type transformations for the generalized hypergeometric function $_{r+2}F_{r+1}(x)$, Zeit. angew. Math. Phys. {\bf 62} (2011) 31--45.

\bibitem{MP2}
A. R. Miller and R. B. Paris, Transformation formulas for the generalized hypergeometric function with integral parameter differences, Rocky Mountain J. Math {\bf 43(1)} (2013) 291--327.

\bibitem{MP3}
A. R. Miller and R. B. Paris, On a result related to transformations and summations of generalized hypergeometric series, Math. Communications {\bf 17} (2012) 205--210.

\bibitem{RR11}
M. A. Rakha and A. K. Rathie,  Generalizations of classical summation theorems for the series ${}_2F_1$ and ${}_3F_2$ with applications, Integral Transforms and Special Functions {\bf 22}, 823--840, 2011.  

\bibitem{RP}
A. K. Rathie and R. B. Paris, Extension of some classical summation theorems for the generalized hypergeometric series with inetegral parameter differences, (2013) submitted for publication.

\bibitem{S}
L. J. Slater, {\it Generalized Hypergeometric Functions}, Cambridge University Press, Cambridge, 1966.

\bibitem{W}
S. Wolfram, The Wolfram Functions Site, http://www.functions.wolfram.com/07.27.17.0034.01 and 07.27.17.0046.01.

}
\end{thebibliography}
\end{document}